\documentclass[a4paper,12pt]{amsart}

\usepackage{url}
\usepackage[margin=0.9in]{geometry}
\usepackage[english]{babel}
\usepackage[utf8]{inputenc}

\usepackage{amsmath}
\usepackage{amssymb}
\usepackage{mathbbol}
\usepackage{algorithm}
\usepackage{algorithmic}
\usepackage{cite}
\usepackage{graphicx}
\usepackage{xcolor}
\usepackage{tikz}
\usepackage{diagbox}
\usepackage{hyperref}

\newtheorem{theorem}{Theorem}[section]

\theoremstyle{definition}
\newtheorem{definition}[theorem]{Definition}
\newtheorem{example}[theorem]{Example}

\theoremstyle{remark}
\newtheorem{remark}[theorem]{Remark}

\numberwithin{equation}{section}

\begin{document}

\title[A Database of Modular Forms on Noncongruence Subgroups]{A Database of Modular Forms on Noncongruence Subgroups}


\author[D. Berghaus]{David Berghaus}
\address{Bethe Center, University of Bonn, Nussallee 12, 53115 Bonn, Germany}
\curraddr{}
\email{berghaus@th.physik.uni-bonn.de}
\thanks{}

\author[H. Monien]{Hartmut Monien}
\address{Bethe Center, University of Bonn, Nussallee 12, 53115 Bonn, Germany}
\curraddr{}
\email{hmonien@uni-bonn.de}
\thanks{}

\author[D. Radchenko]{Danylo Radchenko}
\address{Laboratoire Paul Painlevé, University of Lille, F-59655 Villeneuve d'Ascq, France}
\curraddr{}
\email{danradchenko@gmail.com}
\thanks{}


\date{}

\dedicatory{}

\newcommand{\Hbar}{\overline{\mathcal{H}}}
\newcommand{\bigO}{\mathcal{O}}
\newcommand{\ZZ}{\mathbb{Z}}
\newcommand{\QQ}{\mathbb{Q}}
\newcommand{\RR}{\mathbb{R}}
\newcommand{\CC}{\mathbb{C}}
\newcommand{\subgroup}{\leqslant}
\newcommand{\epsmachine}{\epsilon_\textrm{machine}}
\newcommand{\iu}{{i\mkern1mu}}

\begin{abstract}
	We present a database of several hundred modular forms up to and including weight six on noncongruence subgroups of index $\leq 17$. In addition, our database contains expressions for the Belyi map for genus zero subgroups and equations of the corresponding elliptic curves for genus one subgroups and numerical approximations of noncongruence Eisenstein series to 1500 digits precision.
\end{abstract}

\maketitle
\allowdisplaybreaks

\section{Introduction}
The theory of modular forms on noncongruence subgroups is still not well understood (for recent progress see \cite{MR3788845,Calegari}).
In this work we provide a large number of computed examples in the hope that these will lead to new observations and conjectures. Building a database of modular forms on noncongruence subgroups has been listed as one of the goals for future research in the AIM Workshop \emph{Noncongruence modular forms and modularity} \cite{noncong_workshop}. The usefulness of computer generated data to investigate modular forms on noncongruence subgroups has already been demonstrated by Atkin and Swinnerton-Dyer \cite{asd} who formulated key conjectures based on just a handful of computed examples. Despite this, very few modular forms on noncongruence subgroups have been computed. The known examples were typically restricted to noncongruences character groups \cite{verill,kurth} or to noncongruence subgroups of very low index \cite{fiori_franc} and a few large index genus zero cases \cite{monien_co3}.

In this work we apply the algorithm of \cite{numerics_paper} to make a database of modular forms on noncongruence subgroups. This complements the existing related databases such as the database of classical modular forms on congruence subgroups \cite{lmfdb_modform_db}, the database of Hilbert modular forms \cite{hilbert_form_db}, and the database of Belyi maps \cite{belyi_db}.

\section{Background and Notation}

\subsection{Modular Forms}
Let $\mathrm{SL}(2,\ZZ)$ be the group of $2\times 2$ integer matrices with determinant~$1$ and let $G$ be a finite index subgroup of $\mathrm{SL}(2,\ZZ)$. A modular form of 
weight~$k$ for the group~$G$ is a holomorphic function on the upper half plane $\mathcal{H} := \{\tau \in \CC \, | \, \textrm{Im}(\tau) > 0\}$ that satisfies the transformation law
\begin{equation}
	f(\gamma(\tau)) = (c\tau+d)^k f(\tau)\,,\qquad \mbox{ for all  }\gamma = \begin{pmatrix}
		a & b\\
		c & d
	\end{pmatrix}
	\in G\,,
\end{equation}
where 
\begin{equation}
	\gamma(\tau) := \frac{a\tau +b}{c\tau + d}\,,
\end{equation}
and such that $(c\tau+d)^{-k}f(\gamma(\tau))$ is bounded as $\mathrm{Im}(\tau)\to\infty$ for all $\gamma \in \mathrm{SL}(2,\ZZ)$. In this paper we will restrict to the case when there is no multiplier system and when $k$ is an even integer and without loss of generality we may instead pass to the quotient $\mathrm{PSL}(2,\ZZ):=\mathrm{SL}(2,\ZZ)/\{\pm 1\}$, which we denote by~$\Gamma$ and we consider~$G$ as a finite index subgroup of~$\Gamma$. 
For more insight into the extensive theory of classical modular forms we refer to the books~\cite{cohen_stroemberg2017} by Cohen and Str\"{o}mberg and~\cite{diamond2006first} by Diamond and Shurman, as well as the lecture notes~\cite{Zagier2008} by Zagier. Our setup and notation follows the one used in~\cite{numerics_paper}.

The subgroup $\Gamma(N) \subgroup \Gamma$, defined by
\begin{equation}
	\Gamma(N) := \left\{\begin{pmatrix}
		a & b\\
		c & d
	\end{pmatrix}
	\equiv
	\begin{pmatrix}
		1 & 0\\
		0 & 1
	\end{pmatrix}
	\pmod{N}
	\quad \textrm{and} \quad
	\begin{pmatrix}
		a & b\\
		c & d
	\end{pmatrix}
	\in \Gamma
	\right\}\,,
\end{equation}
is called the \emph{principal congruence subgroup of level N}. A finite index subgroup $G \subgroup \Gamma$ is called \emph{noncongruence} if it does not contain $\Gamma(N)$ for any $N\ge 1$. 

Let $G$ be any finite index subgroup of $\Gamma$ and let $f$ be a modular form of weight $k\in 2\ZZ_{>0}$. Then $f$ has a Fourier expansion
\begin{equation}
	f(\tau) = \sum_{n=0}^{\infty} a_n q_h^n\,,
\end{equation}
where $q_h := \exp(2\pi \iu \tau /h)$ and $h$ denotes the width of the cusp at $\infty$. One can choose bases for spaces of modular forms on noncongruence subgroups with Fourier coefficients defined over $\overline{\QQ}$, and moreover their coefficients can be written in the form (see~\cite{asd,scholl_noncong_modform})
\begin{equation}
	a_n = u^n b_n\,,
\end{equation}
with $b_n$ and $u^h$ belonging to a certain number field $K$. We fix a choice $v$ of a generator for $K$ over $\QQ$ (i.e., $K = \QQ(v)$ with $v$ an algebraic number). We also denote the number field that contains all $a_n$ by $L = \QQ(w)$.

\subsection{Belyi Maps}
\begin{definition}[Belyi Map]
	Let $X$ be a compact Riemann surface. A holomorphic function
	\begin{equation}
		f\, : \, X \rightarrow \mathbb{P}^1(\CC)\,,
	\end{equation}
	is called a \emph{Belyi map} if it is unramified away from three points.
\end{definition}
The covering map
\begin{equation}
	R\, : \, X(G) \rightarrow X(\Gamma) \overset{j}{\cong}
	\mathbb{P}^1(\CC)\,,
\end{equation}
is a Belyi map, where $X(G)=G\backslash \Hbar$ is the modular curve of $G$. If $G$ is a genus zero subgroup then the covering map $R(j_G)$ is a rational function in the Hauptmodul $j_G$ that branches over the images of the elliptic points and the cusps. Note that the Belyi map can be defined over the same field $K$ as above.

\subsection{Signatures and Passports}
\begin{definition}[Signature]
	\label{def:signature}
	We define the \emph{signature} of $G\subgroup \Gamma$ to be the tuple $(\mu,g,n(c),n(e_2),n(e_3))$, where $\mu$ denotes the index of $G$ in $\Gamma$, $g$ denotes the genus of the fundamental domain $\mathcal{F}(G)$, $n(c)$ denotes the number of cusps of $G$ and $n(e_2)$ and $n(e_3)$ denote the number of elliptic points of order two and three, respectively.
\end{definition}
The signature can be used to distinguish between different \emph{types} of subgroups. A further distinction can be made using the notion of \emph{passport}:
\begin{definition}[Passport]
	Given a set of subgroups of equal signature that, by Millington's theorem \cite{millington}, correspond to tuples $(\sigma_S, \sigma_R, \sigma_T)$ which are representatives of $S_\mu$ conjugacy classes, we say that two subgroups belong to the same passport if their \emph{monodromy groups} (i.e., the permutation group generated by $\sigma_S$ and $\sigma_R$) are conjugate in $S_\mu$.
\end{definition}
Passports are useful to know because the action of the absolute Galois group on Belyi maps preserves the passport and hence the size of the passport (i.e., the number of elements inside a passport) provides an upper bound on the degree of the number field $K$. 

\section{Numerical Computation of Elliptic Curves}
\label{sec:elliptic_curve_computation}
In~\cite{numerics_paper} we showed how Fourier expansions of modular
forms on noncongruence subgroups can be computed numerically to high
precision. In this section we show how one can compute from this an
equation for the elliptic curve corresponding to a genus one
noncongruence subgroup. We use an approach similar to that described
by Cremona \cite[Chapter II]{cremona}. Let $f \in S_2(G)$ denote a
nontrivial weight two cusp form for a genus one group~$G$ (note that
for genus one groups $\textrm{dim}(S_2)=1$ which means that~$f$ is
unique up to normalization). The Fourier expansion of~$f$ at all cusps
can be computed numerically by applying the algorithm of
\cite{numerics_paper}. We then define the integral
\begin{equation}
	I_f(\alpha,\beta) := \int_{\alpha}^{\beta} 2\pi\iu f(\tau)d\tau\,,
\end{equation}
where $\alpha,\beta \in \Hbar$ and
\begin{equation}
	I_f(\tau_0) := I_f(\tau_0,\iu\infty) = \int_{\tau_0}^{\iu\infty} 2\pi\iu f(\tau)d\tau = -\sum_{n=1}\frac{a_n}{n}\exp(2\pi\iu n\tau_0/h)\,,
\end{equation}
where $h$ denotes the cusp width. Then the period map $P_f$, defined by
\begin{equation}
	P_f(\gamma) := I_f(\tau_0,\gamma(\tau_0)) = I_f(\tau_0)-I_f(\gamma(\tau_0))\,,
    \qquad \gamma = \begin{pmatrix}
	a & b\\
	c & d
    \end{pmatrix} \in G\,,
\end{equation}
is independent of the choice of the base point $\tau_0$ and defines a homomorphism from $G$ to the period lattice $\Lambda_f = \ZZ w_1 + \ZZ w_2$, with $w_1,w_2 \in \CC$ (see~\cite{cremona}).

\subsection{Evaluation of $P_f(\gamma)$}
\label{sec:P_f_gamma}
Since we work numerically, the base point $\tau_0$ should be chosen in a way that maximizes precision when evaluating $P_f(\gamma)$. This is achieved by maximizing the minimum among the imaginary parts of $\tau_0$ and $\gamma(\tau_0)$. If $G$ has a single cusp, we choose
\begin{equation}
	\tau_0 = -\frac{d}{c} + \frac{1}{c}\iu \in \mathcal{H}\,,
\end{equation}
for $c>0$ (if $c=0$ then $P_f(\gamma) = 0$ and if $c<0$ then we simply choose the other representative $-\gamma$ of $\gamma$ in $\mathrm{PSL}(2,\ZZ)$). Note that for this choice of $\tau_0$ we have $\gamma(\tau_0) = a/c +\iu/c$ which means that the imaginary parts of $\tau_0$ and $\gamma(\tau_0)$ are equal.

If $G$ has multiple cusps, we denote by $A_p$ the cusp normalizer associated to cusp $p$ (using the convention that $A_p(\iu\infty) = p$) and denote cusp expansions at $p$ by $f_p = f|_2 A_p$. Then we obtain
\begin{equation}
	P_{f}(\gamma) = I_{f_p}(\tau_{0,p}) - I_{f_p}((A_p^{-1}\gamma A_p)(\tau_{0,p}))\,.
\end{equation} 
Let $\gamma_p := A_p^{-1}\gamma A_p = \begin{pmatrix}
	a_p & b_p\\
	c_p & d_p
\end{pmatrix}$. Then we choose
\begin{equation}
	\tau_{0,p} = -\frac{d_p}{c_p} + \frac{1}{c_p}\iu \in \mathcal{H}\,,
\end{equation}
and evaluate $P_f(\gamma)$ using the expansion from the cusp $p$ for which the value $c_p\cdot h_p$ is minimized (where $h_p$ denotes the cusp width of cusp $p$).

\subsection{Computation of $\Lambda_f$}
The approach for computing the period lattice $\Lambda_f$ numerically can be summarized as follows:
\begin{enumerate}
	\item Compute the set $S = \{P_f(\gamma_1), P_f(\gamma_2), ...\}$ where $\gamma_j$ denote a set of generators for $G$.
	\item Remove (approximate) zeros (which are trivially on the lattice) and (approximate) duplicate elements from the set~$S$.
	\item Choose any two linearly independent elements $w_1, w_2 \in S$ as a first guess for the basis of $\Lambda_f$.
	\item Express the remaining elements as elements in $\QQ[w_1,w_2]$ by solving
	\begin{equation}
		\begin{pmatrix}
			\textrm{Re}(w_1) & \textrm{Re}(w_2)\\
			\textrm{Im}(w_1) & \textrm{Im}(w_2)
		\end{pmatrix}
		\cdot
		\begin{pmatrix}
			r_j\\
			s_j
		\end{pmatrix}
		=
		\begin{pmatrix}
			\textrm{Re}(P_f(\gamma_j))\\
			\textrm{Im}(P_f(\gamma_j))
		\end{pmatrix}\,,
	\end{equation}
	to obtain $P_f(\gamma_j) = r_j w_1 + s_j w_2$ with $r_j, s_j \in \QQ$.
	\item Find the least common denominator $\lambda$ so that $\lambda r_j, \lambda s_j \in \ZZ$ for all $j$.
	\item Find a basis matrix $\begin{pmatrix}
		a & b\\
		c & d
	\end{pmatrix}$ for the $\ZZ^2$-submodule spanned by $\lambda\cdot (r_j,s_j)$ by computing the Hermite normal form. 
	\item Switch to the lattice basis
	\begin{equation}
		w_1, w_2 \mapsto \frac{a w_1 + b w_2}{\lambda}, \frac{c w_1 + d w_2}{\lambda}\,.
	\end{equation}
	\item Compute $\tau := w_1/w_2$ with $\textrm{Im}(\tau) >0$ (interchanging $w_1$ and $w_2$ if necessary) and transform~$\tau$ to the fundamental domain of $\Gamma$ to obtain a lattice $\ZZ+\ZZ\tau$ equivalent to~$\Lambda_f$.
\end{enumerate}
\begin{remark}
  The precision of the elements
  $S = \{P_f(\gamma_1), P_f(\gamma_2), ...\}$ depends on the choice of
  the generator representatives (which we compute using
  \textsc{Sage}). We are unaware of any algorithm that produces
  generators that are optimal in the sense that $\tau_0$ has maximal
  possible imaginary parts, see Section~\ref{sec:P_f_gamma}. For this
  reason we conjugate $G$ by $\sigma = (1\,j)$, for $j=2,...,\mu$ (in
  terms of the permutation group representation of $G$) and choose the
  representative for which the height of $\tau_0$ becomes
  maximal. While this choice might not be ideal, because it may happen
  that $\tau_0$ is below the height of the fundamental domain, given by $\mathrm{Im}(\mathcal{F}(G)) = \sqrt{3}/(2h_{\mathrm{max}})$ where $h_{\mathrm{max}}$ denotes the largest cusp width (which
  means that we cannot evaluate $P_f(\gamma)$ to full precision), in
  practice the precision has always been sufficient to identify the
  elliptic curve.
\end{remark}

\subsection{Computation of the elliptic curve}
Once the period lattice $\Lambda_f$ has been obtained the computation of the elliptic curve becomes straightforward. We first compute the elliptic invariants
\begin{equation}
	g_2(\tau) = 60 G_4(\tau) = 60 \sum_{(m,n)\neq (0,0)}(m+n\tau)^{-4}\,,
\end{equation}
where $G_4(\tau) = (\pi^4/45) E_4(\tau)$ and
\begin{equation}
	g_3(\tau) = 140 G_6(\tau) = 140 \sum_{(m,n)\neq (0,0)}(m+n\tau)^{-6}\,,
\end{equation}
where $G_6(\tau) = (2\pi^6/945) E_6(\tau)$ using \textsc{Arb} \cite{arb}. Then the $j$-invariant is given by
\begin{equation}
	j(\tau) = 1728\frac{g_2^3}{g_2^3-27g_3^2}\,.
\end{equation}
After that we identify $j(\tau)$ as an algebraic number in $K$ (using
the LLL algorithm \cite{lll}) and use \textsc{Sage} \cite{sage} to
obtain an equation for the elliptic curve.
\begin{example}
  Consider the noncongruence subgroup $G$ with monodromy group
  $\text{PSL}_2\left(\mathbb{F}_8\right)$ and signature $(9, 1, 1, 1, 0)$ that
  corresponds to the admissible permutation pair
  $(\sigma_S, \sigma_R)$ with $\sigma_S= (1)(2\,5)(3\,7)(4\,8)(6\,9)$
  and $\sigma_R= (1\,2\,6)(3\,8\,5)(4\,9\,7)$. We find that
  $\tau = 0.332234...+0.744371...\,\iu$ and the value of the
  $j$-invariant is $\approx -1159088625/2097152$, from which we get
  the elliptic curve
	\begin{equation}
		E \colon y^2 + xy + y = x^3 - x^2 - 95x - 697
	\end{equation}
	of conductor $162$. We can corroborate this result by computing the three weight 6 cusp forms of $S_6(G)$ (which we denote by $x,y,z$ here) and observing that they satisfy the relation
	\begin{equation}
		-xy^2 + (-2w)y^3 + x^2z + (-w)xyz + (19w^2)y^2z + (-3w^3)yz^2 + (15w^4)z^3 = 0\,,
	\end{equation}
	for $w = (-27/256)^{1/9}$. Substituting $x\mapsto xw$, $z\mapsto z/w$ and dividing by $w$ we obtain the cubic
	\begin{equation}
		-xy^2 - 2y^3 + x^2z - xyz + 19y^2z - 3yz^2 + 15z^3 = 0\,,
	\end{equation}
	defined over $\QQ$ from which we obtain the same equation for $E$ using \textsc{Sage}.
\end{example}

\section{Database Structure}
A database of noncongruence subgroups of index $\mu \leq 17$ has been computed by Str\"omberg \cite{stroemberg_recent} and can be obtained at \cite{stroemberg_db}. We have used it as a foundation for our database.
\begin{remark}
    We have developed an algorithm that can list subgroups for indices beyond those computed in \cite{stroemberg_db, stroemberg_recent}. With its help we were able to generate all subgroups of index $\mu \leq 30$. However later we discovered a paper by Vidal \cite{vidal} which describes a much more efficient and essentially optimal algorithm for the same task. We plan to extend the database with higher index subgroups in the foreseeable future.
\end{remark}

\subsection{Database labels}
Each database entry is labeled by the signature of the elements contained in the passport, followed by the id of the passport, as well as a letter that labels the Galois orbit.
\begin{example}
	Consider the passport of signature $(15,1,2,1,0)$ and monodromy group $S_{15}$ which contains three elements (this is the fourth passport of this signature in the database \cite{stroemberg_db} so we give it the id 3). This passport decomposes into two Galois orbits: The first one with label \emph{15\_1\_2\_1\_0\_3\_a} for which $K = \QQ(v)$ where $v^2 - v - 1 = 0$ and the second one with label \emph{15\_1\_2\_1\_0\_3\_b} for which $K = \QQ$.
\end{example}
\subsection{Permutation triple normalization}
Each database entry is unique up to a normalization of the permutation triple (i.e., conjugation in $\textrm{SL}_2(\ZZ)$). We usually normalized~$\sigma_T$ in a way that the cycle type is sorted in decreasing order and that the labels are sorted in increasing order. This means that the largest cusp is placed at infinity. While this causes the number field~$L$ to be of larger degree (and hence the arithmetic to be slower), the factored expressions in~$u$ and~$v$ typically involve smaller factors and are hence are more easily readible.
\begin{example}
	Consider the permutation $(1\,2\,3\,4)(5\,6\,7)(8\,9)$ which has a cycle type of $(4,3,2)$ that is sorted in decreasing order.
\end{example}
There are three exceptions to the above strategy. The first is given for the case when the subgroup~$G$ has multiple cusps of equal width. Placing one of these cusps at infinity may cause~$K$ to be of a larger degree, due to the broken symmetry. In this case we therefore normalize $\sigma_T$ so that the largest cusp with unique cusp width is placed at infinity. The second exception occurs if~$G$ has a cusp of (unique) width~1. In this case we place the cusp of width 1 at infinity because it is convenient that $K = L$. The third exception occurs for genus one subgroups where we do not sort the labels of $\sigma_T$ in increasing order in case another normalization leads to better precision when computing the elliptic curve (see Section~\ref{sec:elliptic_curve_computation}).
\subsection{The data that has been computed for each database element}
For each database element we have computed the curve (either the Belyi map for genus zero subgroups or the elliptic curve for genus one subgroups), as well as the Fourier expansions for the Hauptmoduls, cusp forms, and modular forms up to and including weight~6. For genus zero subgroups we have achieved this by applying the algorithms of~\cite[Section 5]{numerics_paper} over a number field~$L$. For higher genera subgroups we have computed numerical approximations of the corresponding Fourier expansions to 1500 digits of precision using the method of~\cite[Section 4.3.2]{numerics_paper} and from this guessed the closed-form solutions using the LLL algorithm. We remark that for reasons of efficiency, it is not required to compute all forms \emph{from scratch}. Instead, we compute products between cusp forms and modular forms to generate (some of the) higher weight forms. Additionally, we have computed numerical approximations of the basis factors that map $M_k$ onto $E_k$ to 1500 digits precision using an algorithm described in~\cite{eisenstein_paper}. Interestingly, we were unable to find any algebraic dependencies among these expressions.
\subsection{A complete example}
For an explicit example, let us take a look at the database entry $15\_1\_2\_1\_0\_2\_a$.
\begin{itemize}
	\item $G$: The subgroup $G$ corresponds to the permutation group generated by 
    \begin{align*}
    \sigma_S &= (1\,15)(2\,12)(3\,7)(4\,9)(5\,13)(6\,10)(8\,14)(11)\,,\\ 
    \sigma_R &= (1\,11\,12)(2\,13\,6)(3\,8\,15)(4\,10\,7)(5\,14\,9)\,,\\
    \sigma_T &= (1\,3\,4\,5\,6\,7\,8\,9\,10\,2)(11\,12\,13\,14\,15)\,.
    \end{align*}
	\item Monodromy Group: $(C_3 \times C_3 \times C_3 \times C_3) : (C_2 \times A_5)$.
	\item $K$: $v=1$ (i.e., $K=\QQ$).
	\item Embeddings: $v=1$.
	\item $u$: $(-3125/14348907)^{1/10}$ with an embedding $-0.409269... - 0.132979...\cdot\iu$.
	\item $L$: $w^8 - 5/27w^6 + 25/729w^4 - 125/19683w^2 + 625/531441$ with $w = -0.409269... - 0.132979...\cdot\iu$.
	\item Curve: $y^2 + xy + y = x^3 - x^2 + 7x - 103$.
	\item Fourier expansions (up to $\sim 700$ terms):
	\begin{itemize}
		\item $M_2$: 
		\begin{itemize}
			\item $m_0 = 1 + 48/5u^2q_{10}^2 - 528/5u^3q_{10}^3 - 2448/25u^4q_{10}^4 - 3096/25u^5q_{10}^5 + ...$.
			\item $m_1 = q_{10} + uq_{10}^2 - 21/5u^2q_{10}^3 + 101/5u^3q_{10}^4 + 816/25u^4q_{10}^5 + ...$.
		\end{itemize}
		\item $S_2$:
		\begin{itemize}
			\item $s_0 = q_{10} + uq_{10}^2 - 21/5u^2q_{10}^3 + 101/5u^3q_{10}^4 + 816/25u^4q_{10}^5 + ...$ (this is the same as $m_1$ in this case).
		\end{itemize}
		\item $E_2$:
		\begin{itemize}
			\item $e_0$: We give the Eisenstein basis matrix $[1,0.356085...+ 0.115699...\cdot \iu]$ (up to 1500 terms), meaning that $e_0 = m_0 + (0.356085...+ 0.115699...\cdot \iu)m_1$. Note that this basis is already in canonical normalization.
		\end{itemize}
		\item $M_4$:
		\begin{itemize}
			\item $m_0 = 1 - 688747536/625u^{10}q_{10}^{10} + ...$ (this is an oldform from $\Gamma$, namely the weight 4 Eisenstein series for $\mathrm{SL}(2,\ZZ)$).
			\item $m_1 = q_{10} - 70111/110u^4q_{10}^5 + ...$.
			\item $m_2 = q_{10}^2 - 621/22u^3q_{10}^5 + ...$.
			\item $m_3 = q_{10}^3 + 237/22u^2q_{10}^5 + ...$.
			\item $m_4 = q_{10}^4 - 115/22uq_{10}^5 + ...$.
		\end{itemize}
		\item $S_4$:
		\begin{itemize}
			\item $s_0 = q_{10} - 65u^3q_{10}^4 - 1488/5u^4q_{10}^5 + ...$.
			\item $s_1 = q_{10}^2 - 27/5u^2q_{10}^4 + ...$.
			\item $s_2 = q_{10}^3 - uq_{10}^4 + 16u^2q_{10}^5 + ...$.
		\end{itemize}
		\item $E_4$:
		\begin{itemize}
			\item $e_0 = m_0$.
			\item $e_1 = m_1 + (12.068314...+3.921233...\cdot\iu)m_2 + (106.205344...+77.162699...\cdot\iu)m_3 + (-61.318023...-84.397019...\cdot\iu)m_4$.
			\item $e_{0,\textrm{can}}$: For the canonical normalizations we get\\
			$e_{0,\textrm{can}} = e_0 + (-0.014977...-0.004866...\cdot\iu)e_1$.
			\item $e_{1,\textrm{can}} = (0.059909...+ 0.019465...\cdot\iu)e_1$.
		\end{itemize}
	\end{itemize}
\end{itemize}
(The weight six spaces are also given in the database, but are not listed here.)

\subsection{Status of the Database}
The number of computed passports in the current version can be found in Table~\ref{tab:amount_of_computed_pps}.
\begin{table}
	\centering
	\begin{tabular}{|l||*{2}{c|}}\hline
		\backslashbox{index}{genus}
		&\makebox[3em]{0}&\makebox[3em]{1}\\\hline\hline
		7 & 3/3 & 0\\\hline
		8 & 1/1 & 0 \\\hline
		9 & 9/9 & 1/1 \\\hline
		10 & 9/9 & 1/1 \\\hline
		11 & 6/6 & 0 \\\hline
		12 & 23/27 & 2/3 \\\hline
		13 & 22/23 & 1/1 \\\hline
		14 & 21/29 & 1/2 \\\hline
		15 & 54/62 & 7/9 \\\hline
		16 & 36/65 & 7/9 \\\hline
		17 & 16/35 & 1/2 \\\hline
		\hline
		total & 200/269 & 21/28\\
		\hline
	\end{tabular}
	\caption{Number of computed passports that are currently in the database. (Note that there are no noncongruence subgroups with $\mu \leq 17$ and $g > 1$.)}
	\label{tab:amount_of_computed_pps}
\end{table}
At the current stage the size of the database is $\sim 6 \textrm{GB}$ in compressed and $\sim 16 \textrm{GB}$ in uncompressed form. In total, about $25\,000$ hours of CPU time on \texttt{Intel Xeon E5-2680 v4 @ 2.40GHz} was used to compute the data.

\subsection{Reliability of the Results}
For genus zero, all coefficients have been computed using rigorous arithmetic over number fields or rigorous interval arithmetic. The only exception are the Eisenstein series whose precision we can only estimate heuristically. For higher genera, the coefficients are non-rigorous (but supported by highly convincing numerical evidence) because they have been guessed using the LLL algorithm. As an additional verification, we have also compared the numerical values of the rigorous expressions and the Eisenstein series to a computation using Hejhal's method with a different horocycle height to verify that the results match to at least 10 digits of precision.

\subsection{How to access the Database}
The database is currently available as a \textsc{GitHub} repository at \cite{noncong_database_github} and planned to be released to the LMFDB soon. The database entries can be loaded by running \textsc{Sage} scripts which return a \textsc{Python} dictionary containing the results. This results in more portability between different versions than storing the results as pickled objects. Additionally we provide the results in printed form as strings inside \textsc{JSON} files. To save on storage, we did not store the numerical approximations to the Fourier coefficients. For the same reason we also did not store the expressions over the number field $L$ explicitly but instead generate them when loading the \textsc{Sage} scripts by plugging in the values of~$v$ and~$w$.

\section{Some notable examples}
\subsection{The largest degree of $K$}
The largest degree of $K$ occurs for the passport $16\_0\_3\_2\_1\_11\_a$ for which we get $K = \QQ(v)$ where $v=1.531805...+0.185685...\iu$ is a root of
\begin{align*}
	0 &= v^{24} - 6v^{23} + 21v^{22} - 60v^{21} + 184v^{20} - 478v^{19} + 651v^{18} - 1220v^{17} + 2230v^{16} + 1226v^{15} \\
	 &- 947v^{14} + 804v^{13} - 7092v^{12} - 6862v^{11} + 3971v^{10} - 15340v^9 + 7975v^8 + 36044v^7 + 7134v^6 \\
	 &+ 14896v^5 + 13928v^4 - 2372v^3 + 3970v^2 - 584v + 22\,,
\end{align*}
with discriminant $2^{52} 3^{15} 5^{10} 7^4 13^4$ and Galois group $S_{24}$.
While it is definitely possible to use the method of \cite[Section 5.1]{numerics_paper} to compute Belyi maps for higher number field degrees (see for example \cite{monien_co3}), computing non-trivial amounts of Fourier coefficients becomes infeasible, which is the reason why these large passports have not been included in the database. To give an idea of the size of the coefficients involved, we remark that the first non-trivial coefficient of the Hauptmodul starts with
\begin{align*}
	18035144800333385601165709955931694753948569033237338483618065475555557037...\\
	...3958967421781810036632024163158272374698352311806451683704v^{23}+...\,,
\end{align*}
which has 132 decimal digits and factors into 
\begin{align*}
2^3 \cdot 3^4 \cdot 7^2 \cdot 13^3 \cdot 179 \cdot 1178062360621513 \\ \cdot\,1515687995725535658492175921 \cdot 15731038963473855359968899262003 \\ \cdot\,514197500928774955284304452843050396002488377491\,.
\end{align*}
\begin{figure}
	\centering
	\includegraphics[width=\columnwidth]{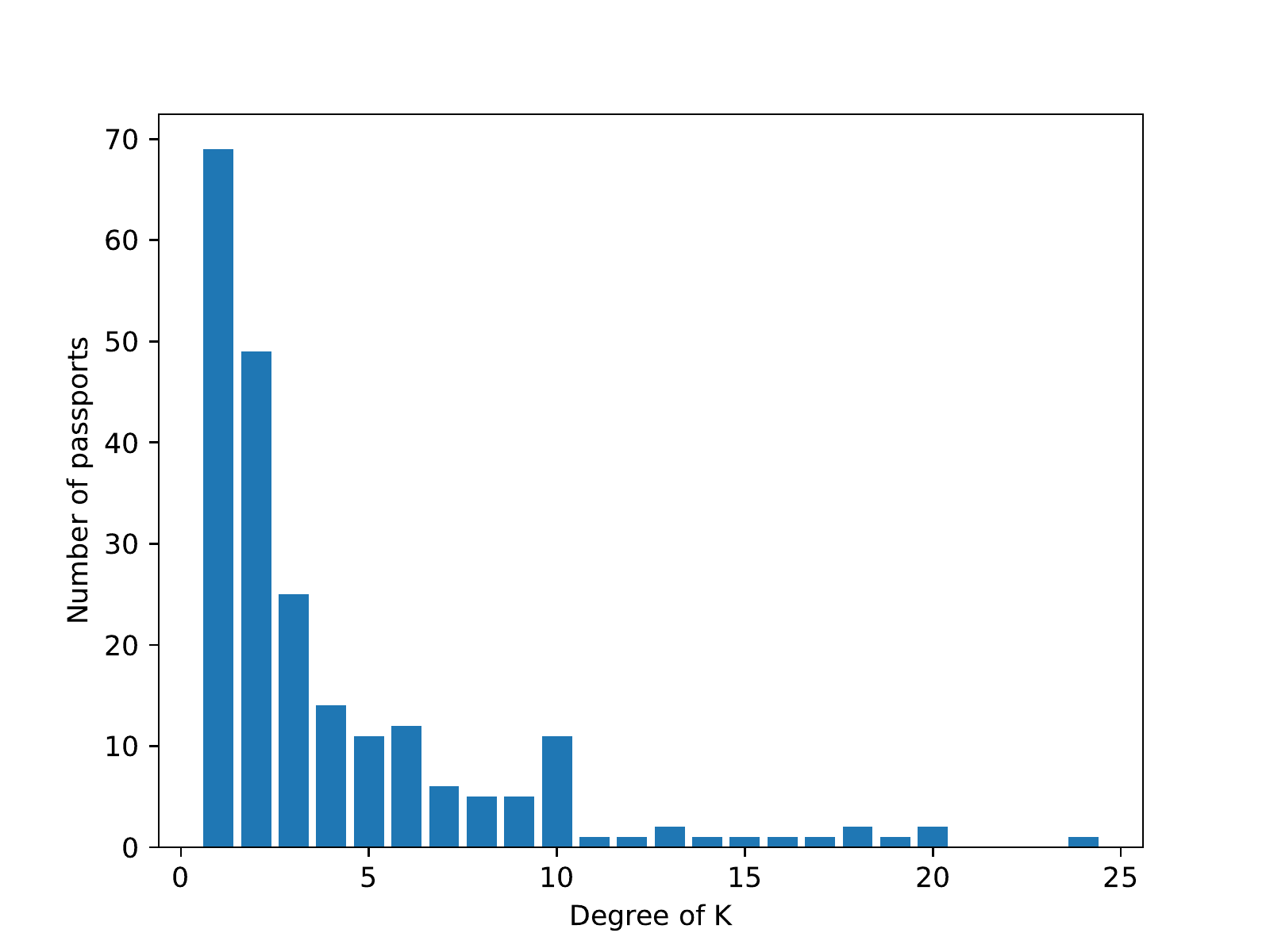}
	\caption{Distribution of the degrees of $K$ in the database.}
	\label{fig:K_degree_histo}
\end{figure}
The distribution of degrees of $K$ is given in Figure~\ref{fig:K_degree_histo}.
\subsection{Most Fourier expansion terms}
\begin{figure}
	\centering
	\includegraphics[width=\columnwidth]{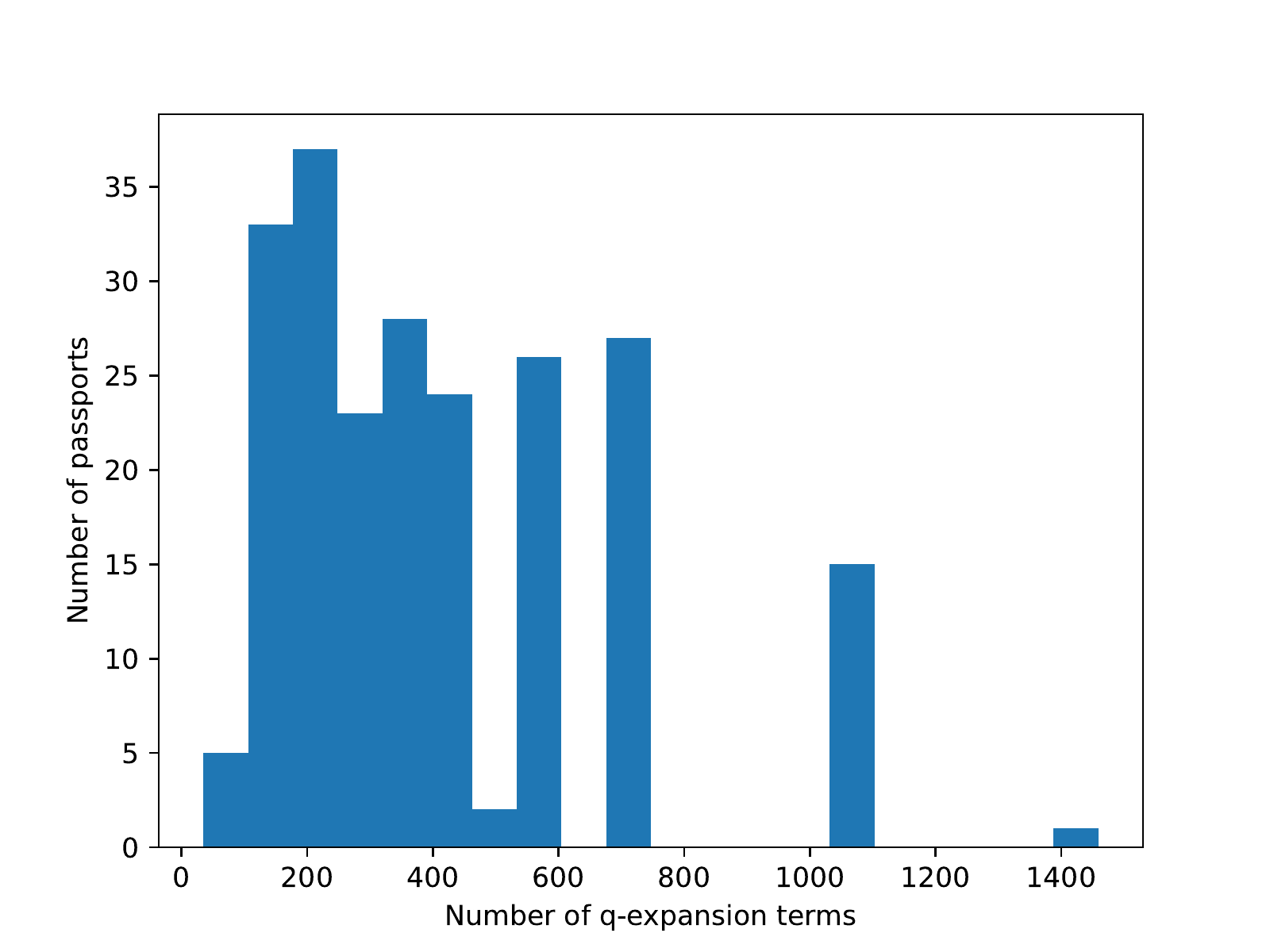}
	\caption{Distribution of the number of terms of the Fourier expansions in the database.}
	\label{fig:precs_histo}
\end{figure}
The largest number of Fourier expansion coefficients that have been computed for the database occurs for the passport $9\_1\_1\_1\_0\_0\_a$ for which we list 1459 coefficients. 
A distribution of the number of Fourier expansion terms is given in Figure~\ref{fig:precs_histo}. For genus zero subgroups we chose the number of terms based on some heuristic guesses, which depend on the degree of $L$, ensuring that the computation does not last longer than a few days. For higher genera subgroups we included as many terms as the LLL could determine reliably from the numerical approximations. For certain passports we also had to reduce the number of coefficients in order to satisfy the maximal \textsc{GitHub} file size limit of 100MB.

\subsection{Elliptic curves defined over $\QQ$}
\begin{table}
	\centering
	\begin{tabular}{c|c|c}
		Passport Label & Elliptic Curve & Conductor\\
		\hline
		$9\_1\_1\_1\_0\_0\_a$ & $y^2 + xy + y = x^3 - x^2 - 95x - 697$ & 162\\
		$10\_1\_1\_0\_1\_0\_a$ & $y^2 + xy + y = x^3 + x^2 - 110x - 880$ & 15\\
		$15\_1\_2\_1\_0\_2\_a$ & $y^2 + xy + y = x^3 - x^2 + 7x - 103$ & 270\\
		$15\_1\_2\_1\_0\_3\_b$ & $y^2 + xy = x^3 - x^2 - 41370x + 2022196$ & 2970
	\end{tabular}
	\caption{Noncongruence subgroups corresponding to elliptic curves over $\QQ$.}
	\label{tab:QQ_elliptic_curves}
\end{table}
In Table~\ref{tab:QQ_elliptic_curves} we give all of the examples of genus~1 noncongruence subgroups from the database that correspond to elliptic curves over~$\QQ$, together with their defining equations and their conductors.


\bibliographystyle{abbrv}
\bibliography{references.bib}
\end{document}